\documentclass[12pt]{article}
\usepackage{stmaryrd}
\usepackage{amssymb}
\usepackage{color, amsmath,amssymb, amsfonts, amstext,amsthm, latexsym}

\usepackage{amssymb, epsfig, amssymb, latexsym}
\usepackage{amsmath}
\usepackage{graphicx}
\usepackage{longtable}
\usepackage{graphicx}
\usepackage{subfigure}
\usepackage{appendix}
\usepackage{cite}
\usepackage[colorlinks=false, linktocpage=true]{hyperref}
\usepackage{authblk}

\numberwithin{equation}{section} \topmargin -0.4in

\renewcommand{\phi}{\varphi}

\title{State transitions in the Morris-Lecar model under stable L\'evy noise
  }

\author[1,3]{Rui Cai}
\author[2,3]{Yancai Liu \thanks{Corresponding author: yancailiu@hust.edu.cn}}
\author[4]{Jinqiao Duan}
\author[2,3]{Almaz Tesfay Abebe}
\affil[1]{School of Science, Hubei University of Technology, Wuhan 430068, China}
\affil[2]{School of Mathematics and Statistics, Huazhong University of Science and Technology, Wuhan 430074, China}
\affil[3]{Center for Mathematical Sciences, Huazhong University of Science and Technology, Wuhan 430074, China}
\affil[4]{Department of Applied Mathematics, Illinois Institute of Technology, Chicago, IL 60616, USA}

\date{}

\begin{document}
\maketitle

\begin{abstract}
This paper considers the state transition of the stochastic Morris-Lecar neuronal model driven by symmetric $\alpha$-stable L\'evy noise. The considered system is bistable: a stable fixed point (resting state) and a stable limit cycle (oscillating state), and there is an unstable limit cycle (borderline state) between them. Small disturbances may cause a transition between the two stable states, thus a deterministic quantity, namely the maximal likely trajectory, is used to analyze the transition phenomena in non-Gaussian stochastic environment. According to the numerical experiment, we find that smaller jumps of the L\'evy motion and smaller noise intensity can promote such transition from the sustained oscillating state to the resting state. It also can be seen that larger jumps of the L\'evy motion and higher noise intensity are conducive for the transition from the borderline state to the sustained oscillating state. As a comparison, Brownian motion is also taken into account. The results show that whether it is the oscillating state or the borderline state, the system disturbed by Brownian motion will be transferred to the resting state under the selected noise intensity.
\end{abstract}


\section{Introduction}

Recent studies have shown that neuroelectrophysiological activities contain complex nonlinear dynamic behaviors \cite{izhikevich2007dynamical}. It is impossible to fully explain the various phenomena in neuroscience by simply using traditional linear viewpoints and statistical methods, and the simple description of the experimental results can not meet the requirements of quantitative analysis for neuroscience. Therefore, based on a large amount of experimental data, it is important to apply nonlinear dynamics and mathematical methods to model and analyze neuronal systems, which promotes the development of computational neuroscience. Computational neuroscience is established on the basis of biology and build mathematical model based on biological knowledge, and simulation results are obtained through numerical simulation, so as to study biological physiological characteristics. The advances in computational neuroscience could bring humans closer to preventing or alleviating diseases like Parkinson's and depression \cite{Chakravarthy2009Homeostasis}.

%As we all know, computational neuroscience is an interdisciplinary aspect which is based on biological models and simulates nervous activities through various numerical calculation methods to reveal the biological physiological characteristics.

In the study of neuroscience, it is a remarkable fact that neurons live in a noisy environment and are affected by various noises \cite{ermentrout2010mathematical}. It is inevitable to consider the interference of noise when studying some behaviors of neurons. The sources of noise include random opening and closing of ion channels, and depolarization and hyperpolariations bursts caused by spontaneous release of neurotransmitters. In medicine, research on the treatment of diseases such as epilepsy and Parkinson by applying external electric field stimulation to the human body is increasing. Therefore, it is necessary to take neuron model as the research object, and accurate mathematical model obtained by adding noise terms is the key to the study of neuronal system problems. There is a lot of research on the neural model under the influence of noise, mostly focusing on Gaussian white noise \cite{franovic2015activation,FFranovic2015activation,wang2017electrical,longtin2013neuronal,lindner2004effects}. For the non-Gaussian L\'evy noise which is more extensive than Gaussian noise, developments in the related stochastic neuronal system\cite{patel2008stochastic,patel2007levy,wang2016levy} are relatively small. It is noteworthy that recent empirical research has shown that the probability distribution of ‘anomalous’ (high amplitude) neural oscillations has heavier tail than the standard normal distribution\cite{roberts2015heavy}. Therefore, L\'evy motion is well suited to modelling such kind of noise.

The existing research shows that noise can cause various phenomena in dynamical systems \cite{luchinsky1999observation}, such as stochastic resonance \cite{perc2007stochastic,ishimura2016stochastic,duki2018stochastic}, chaos \cite{bashkirtseva2013noise,frey1993noise,bashkirtseva2016stochastic} and state transitions\cite{xu2013levy,zheng2016transitions,Wu2018L}. Especially, the phenomena of noise induced state transitions in neuronal systems are widely found \cite{lim2010noise,tanabe2001noise,touboul2012mean}. So we focus our attention in this paper on the state transition of the stochastic Morris-Lecar neuron model described by a dynamical system driven by the symmetric $\alpha$-stable L\'evy noise. The Morris-Lecar model is a two-dimensional biological neuron model used to reproduce various oscillatory behaviors associated with $Ca^{2+}$ and $K^{+}$ conductance in giant barnacle muscle fibers \cite{morris1981voltage}. This model is simple in form, low in order, and comprehensively reflects the various characteristics of neurons. It plays an increasingly important role in the field of neuroscience. So much attention has been payed on deterministic or stochastic Morris-Lecar model \cite{Terman1992The,newby2013breakdown,tateno2004random}. In \cite{morris1981voltage}, the authors showed that the stochastic Morris-Lecar neuron can be approximated by a two-dimensional Ornstein-Uhlenbeck (OU) modulation of a constant circular motion in a neighborhood of its stable point. In \cite{keener2011perturbation}, a stochastic interpretation of spontaneous action potential initiation was developed for the Morris-Lecar equations. However, related research on the effects of non-Gaussian L\'evy noise on the Morris-Lecar model is relatively rare.

The deterministic Morris-Lecar model is represented by the following second-order system:
\begin{align}\label{DE}
&d{v_{t}} = \frac{1}{C}[-g_{Ca}m_{\infty}(v_{t})(v_{t}-V_{Ca})-g_{K}w_{t}(v_{t}-V_{K})-g_{L}(v_{t}-V_{L})+I]dt,\\
&d{w_{t}} =\phi\frac{w_{\infty}(v_{t})-w_{t}}{\tau_{w}(v_{t})}dt,
\end{align}
where
\begin{align}
&m_{\infty}(v)=0.5[1+\tanh(\frac{v-V_{1}}{V_{2}})],\\
&w_{\infty}(v)=0.5[1+\tanh(\frac{v-V_{3}}{V_{4}})],\\
&\tau_{w}(v)=[\cosh(\frac{v-V_{3}}{2V_{4}})]^{-1}.
\end{align}
Here, the variable $v_{t}$ represents the membrane potential, $w_{t}$ is a recovery variable, which represents the evolution of the potassium ion channel open probability. The parameter $C$ is the membrane capacitance, $\phi$ represents the change between the fast and slow scales of neurons. The parameters $g_{Ca},~g_{K},~g_{L}$ are the maximum conductance of calcium, potassium and leakage current channels, respectively. And $V_{Ca},~V_{K},~V_{L}$ are the reversal potential of the above channels, respectively. The parameter $I$ represents the total synaptic input current from the environment, $m_{\infty}$ and $w_{\infty}$ are the steady-state values of the opening probability of $Ca^{2+}$ channel and $K^{+}$ channel, respectively. The parameters $V_{1},V_{2},V_{3},$ and $V_{4}$ are parameters chosen to fit voltage-clamp data.

\begin{figure}[ht]
\centering %使插入的图片居中显示
\includegraphics[height=6cm ,width=8cm]{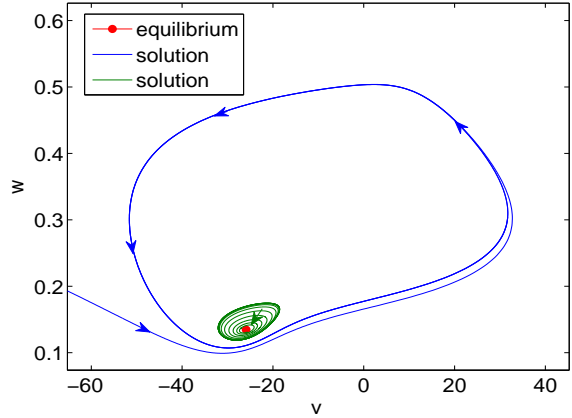}
\caption{Dynamical landscape of the Morris-Lecar model. One stable state (red point) is a spiral sink which coexists with a limit cycle, the other stable state is the big stable limit cycle (blue). The unstable limit cycle (green curve) separates those initial conditions that approach the stable fixed point from those that approach the big stable limit cycle. The arrow indicates the direction of the limit cycle. The parameter values for the type \uppercase\expandafter{\romannumeral2} excitability of ML model are: $C = 20~\mu F/cm^2$, $V_{Ca} = 120~mV$, $V_K =-84~mV$, $V_L = -60~mV$, $g_{Ca} = 4.4~\mu S/cm^2$, $g_K = 8~\mu S/cm^2$, $g_L = 2~\mu S/cm^2$, $V_1 = -1.2~mV$, $V_2 = 18~mV$, $V_3 = 2~mV$, $V_4 = 30~mV$ and $\phi = 0.04$, $I = 92~\mu A/cm^2
$.} % 插入图片的标题，一般放在图片的下方，放在表格的上方
\label{fig:1}
\end{figure}

Figure \ref{fig:1} shows the phase portrait of the deterministic Morris-Lecar model, where the parameters are shown below the figure. It follows from the Figure \ref{fig:1} that the Morris-Lecar model is bistable, there exist both a unique stable fixed point and a stable limit cycle (the big limit cycle in Figure \ref{fig:1}). Also, it is worth noting that there also exists an unstable periodic solution (the small limit cycle in Figure \ref{fig:1}). This trajectory separates those initial conditions that approach the stable fixed point from those that approach the big stable limit cycle. The stable limit cycle corresponds to the sustained oscillating state of neuronal system described by the Morris-Lecar model, while the unique fixed point corresponds to the resting state of the neuron.

In this paper, our main consideration is whether the Morris-Lecar model will have state transition between oscillating state and resting state, under the influence of non-Gaussian L\'evy noise. This phenomenon is considered from two aspects:

1) Whether the solution trajectory starting from the stable limit cycle will enter the attraction basin (ie, the green area) of the stable point under the interference of L\'evy noise. In other words, does the L\'evy noise cause the Morris-Lecar system to shift from the sustained oscillating state to the resting state?

2) Which attraction basin does the solution trajectory starting from the unstable periodic solution (borderline state) will enter under the interference of L\'evy noise? The sustained oscillating state or the resting state? That is to say which state will the borderline state under the influence of non-Gaussian noise become?

Dynamical system theory provides a powerful tool for analyzing deterministic systems of nonlinear differential equations, including related models of neuroscience \cite{Liu2014,Tsumoto2006}. In these theories, the solution is regarded as a curve in phase space, and the dynamic behavior of the model is studied by exploring various geometric structures of the solution. To investigate the behavior of the trajectories in stochastic environment, an indicator should be introduced: the maximal likely trajectory. Then we will use it to explore the state transition behavior of the Morris-Lecar model under the interference of L\'evy noise as well as Brownain noise.

This paper is organized as follows: in Section \ref{TSM}, the stochastic Morris-Lecar model driven by symmetric $\alpha$-stable L\'evy noise is presented. In Section \ref{TM}, we present the method that calculats the maximal likely trajectory and the probability density function diagram. In Section \ref{R}, the effects of L\'evy motion index and noise intensity on the maximal likely trajectories have been displayed. Finally, main conclusions are summarized in Section \ref{CON}.

\section{The Stochastic Model}\label{TSM}

Our purpose in this work is to analyze the state transition of the stochastic Morris-Lecar model. Inspired by \cite{tateno2004random}, the perturbations are represented by a non-Gaussian L\'evy noise current added to the membrane voltage. The dynamics of the Morris-Lecar model driven by such a stimulation is described by the following stochastic differential equations:
\begin{equation}\label{SV}
\begin{split}
&d{v_{t}} =\frac{1}{C}[-g_{Ca}m_{\infty}(v_{t})(v_{t}-V_{Ca})-g_{K}w_{t}(v_{t}-V_{K})-g_{L}(v_{t}-V_{L})+I]dt+\sigma dL_{t}^{\alpha},\\
&d{w_{t}} =\phi\frac{w_{\infty}(v_{t})-w_{t}}{\tau_{w}(v_{t})}dt,
\end{split}
\end{equation}
where $L_{t}^{\alpha}$ is symmetric $\alpha$-stable L\'evy motions, and $\sigma$ is nonnegative noise intensity. The L\'evy motion, as an appropriate model for non-Gaussian processes with jumps \cite{sato1999levy,bertoin1998levy}, has the properties of stationary and independent increments as Brownian motion. It is worth mentioning that the distribution for a stable random variable is denoted as $S_{\alpha}(\delta,\beta,\gamma)$, here $\alpha$ is called the L\'evy motion index (non-Gaussianity index), $\delta$ is the scale parameter, $\beta$ is the skewness, and $\gamma$ is the shift. The following is the definition of a symmetric $\alpha$-stable L\'evy motion.

A $\alpha$-stable L\'evy motion $L_t^\alpha$, with $0 < \alpha  < 2$,  is a stochastic process with the following properties \cite{applebaum2009levy,duan2015introduction}:

(i) $L_0^\alpha$=0, almost surely  (a.s);

(ii) $L_t^\alpha$ has independent increments;

(iii) $L_t^\alpha-L_s^\alpha \sim S_{\alpha}((t-s)^{\frac{1}{\alpha}}, 0, 0)$;

(iv) $L_t^\alpha$ has stochastically continuous sample paths: for every $s$,   $L_t^\alpha \to L_s^\alpha$ in probability, as $t \to  s$.

When $\alpha$ is close to 0, $\alpha$-stable L\'evy motion has larger jumps with lower jump probabilities, while it has smaller jumps with higher jump frequencies for $1<\alpha<2$. Moreover, the jump measure is defined as follows \cite{samorodnitsky1996stable,sato1999levy}:
\begin{align}
  \nu_{\alpha}=\frac{C_{\alpha}dy}{|y|^{1+\alpha}}.
\end{align}
For $0<\alpha<2$, the following tail estimate holds \cite{samorodnitsky1996stable}:
\begin{align}
  \lim_{y\rightarrow\infty}y^{\alpha}\mathbb{P}(L_{t}^{\alpha}>y)=C_{\alpha}\frac{1+\beta}{2}\sigma^{\alpha},
\end{align}
where $C_{\alpha}$ is a positive constant. Unlike Brownian motion's tail estimate which decays exponentially, this estimate indicates that L\'evy motion $L_{t}^{\alpha}$ has a ``heavy tail" which decays polynomially.

\section{The method}\label{TM}

In deterministic dynamical systems, the relevant dynamical behavior of the model can be predicted by studying the geometry of the solution in the phase space. While, for stochastic dynamic systems, the maximal likely trajectory as a connecting tool extends the phase portrait concept of deterministic dynamical systems to stochastic dynamical systems. Therefore, it can help us explore some dynamic behaviors of the Morris-Lecar system under the influence of non-Gaussian L\'evy noise.

Using the method of maximal likely trajectory to explore the impact of L\'evy noise on the stochastic dynamical system has achieved some results. The stochastic pitchfork bifurcation for a system under multiplicative stable L\'evy noise by the maximal likely trajectory is studied in \cite{wang2018a}.  Similarly, this indicator is also used in gene regulatory systems under stable L\'evy noise \cite{chen2019most}. The definition of the maximal likely trajectory is as follows:

The maximal likely trajectory \cite{duan2015introduction}: we consider the solution $S_{t}$ of a stochastic dynamical system starting at the initial points $s_{0}$ in the state space $R^{n}$. At a given time instant $t$, the maximizer $s_{m}(t)$ for $p(s,t)$ which is the probability density function of the solution $S_{t}$ indicates the most probable location of this trajectory at time $t$. The trajectory traced out by $s_{m}(t)$ is called the maximal likely trajectory starting at $s_{0}$. In our case, the state space $R^{n}$ is $R^{2}$.

In order to describe the state transition in the stochastic Morris-Lecar system by the maximal likely trajectory, the first step is to solve the corresponding nonlocal Fokker-Planck equation. Generally speaking, like (\ref{SV}), the two-dimension stochastic differential equation has the following form:
\begin{eqnarray}
% \nonumber to remove numbering (before each equation)
  dv_{t} &=&f_{1}(v_{t},w_{t})dt+\sigma dL_{t}^{\alpha}, \nonumber \\
  dw_{t} &=&f_{2}(v_{t},w_{t})dt. \nonumber
\end{eqnarray}
The Fokker-Planck equation of this system is:
\begin{equation}\label{PRO}
p_{t}(v,w,t)=-(f_{1}p)_{v}-(f_{2}p)_{w}+
\sigma^{\alpha}\int_{\mathbb{R}\setminus{\{0\}}}[p(v+v',w,t)-p(v,w,t)]\nu_{\alpha}(dv').\\
\end{equation}
with the initial condition：$p(v,w,0)=\delta(v-v_{0},w-w_{0})$.

In this paper, we adopt the finite difference method proposed in \cite{GAO20161} to discretize the Fokker-Planck equation for a numerical simulation. In the following, we present the numerical algorithms the case of $D = (a; b)\times(c; d)$ with the natural boundary condition. In this case, we can simplify equation (\ref{PRO}) as:
\begin{align}
 p_t(v,w,t) =&-(f_1p)_v - (f_2p)_w  +\sigma_{1}^{\alpha}C_\alpha \int_{a-v}^{b-w} \frac{p(v\!+\!v', w, t) \!-\! p(v, w, t)}{|v'|^{1\!+\!\alpha}}\; {\rm d}v', \nonumber
\end{align}
for ~$(v,w)\in D$; and $p(v,w,t)=0$ for ~$(v, w) \notin D$.

Set ~$x=\frac{2}{b\!-\!a}(v\!-\!a)\!-\!1$,~$y=\frac{2}{d\!-\!c}(w\!-\!c)\!-\!1$, $p(\frac{b\!-\!a}{2}x\!-\!\frac{a\!-\!b}{2},
\frac{d\!-\!c}{2}y-\!\frac{c\!-\!d}{2},t)=P(x,y,t)$,~i.e.,~$v=\frac{b\!-\!a}{2}x\!-\!\frac{a\!-\!b}{2}$, $w=\frac{d\!-\!c}{2}y\!-\!\frac{c\!-\!d}{2}$, then we get
\begin{align}\label{FF}
P_t(x,y,t) =&-\frac{2}{b-a}(f_1 P)_x  -\frac{2}{d-c} (f_2 P)_y \nonumber  \\
 &+ C_\alpha (\frac{2\sigma_1}{b\!-\!a})^{\alpha}\int_{-1\!-\!x}^{1\!-\!x} \frac{P(x\!+\!x',s,t) \!-\! P(x,y,t)}{|x'|^{1\!+\!\alpha}}\; {\rm d}x',
\end{align}
for ~$(x, y) \in D'=(-1,1)\times(-1,1)$; $P(x,y,t)=0$,~$(x,y) \notin  D'$.

Then, we use a numerical method to discretize the nonlocal Fokker-Planck equation (\ref{FF}). We divide the interval~$[-2,2]\times[-2,2]$~in space into $(4J)^{2}$ sub-intervals and define~$x_i=ih$, $y_j=jh$ for $-2J\leq i,j \leq 2J$,~$h=1/J$. Denoting the numerical solution of $P$ at $(x_i, y_j, t)$ by $P_i,j$, we obtain the semi-discrete equation:
\begin{align}
\dfrac{dP_{i,j}}{dt} =& C_{hx} \frac{P_{i-1,j} - 2P_{i,j} + P_{i+1,j}}{h^2}
-\frac{2}{b-a}[(f_1P)_{x,i}^+ + (f_1P)_{x,i}^-]-\frac{2}{d-c}[(f_2P)_{y,j}^+ +  (f_2P)_{y,j}^-] \nonumber  \\
&+ (\frac{2\sigma_{1}}{b-a})^{\alpha} C_{\alpha} h \sum_{k_1=-J-i,k_1 \neq 0}^{J-i}{''}
{\frac{P_{i+k_1,j} - P_{i,j}}{|{x_{k_{1}}}|^{1+\alpha}} }
\end{align}
where $0<\alpha<2$, $i,j = -J+1, \cdots, -2,-1,0,1,2, \cdots, J-1$, ~$C_{hx}\!= -(\frac{2\sigma_1}{b-a})^{\alpha} C_{\alpha}\zeta(\alpha-1) h^{2-\alpha}$~ are constants, $\zeta$ is the Riemann zeta function, the superscripts $\pm$ denote the global Lax-Friedrichs flux splitting defined as $(f_kP)^{\pm}_{i,j} =
\frac{1}{2}(f_kP_{i,j} \pm \alpha_k P_{i,j})$ with $\alpha_k=\max |f_k(x,y)|$,~$k=1,2$. The
summation symbol $\sum{''}$ means that the quantities corresponding
to the two end summation indices are multiplied by $1/2$.

When $L_{t}^{\alpha}$ is replaced by Brownian motion $B_{t}$, the equation becomes
\begin{equation}
p_{t}(v,w,t)=-(f_{1}p)_{v}-(f_{2}p)_{w}+\frac{\sigma^{2}}{2}p_{vv}.\\
\end{equation}

For the stochastic differential equations driven by L\'evy process, the corresponding Fokker-Planck equation includes a non-local term, namely the fractional Laplacian term, which quanifies the effects of non-Gaussian. Since the equation (\ref{PRO}) satisfies the condition of the theorem in \cite{wei2015well}, the equation (\ref{PRO}) under consideration has a weak $L^{p}$ solution. It is known that the analytic solutions for this kind of Fokker-Planck equations are difficult to obtain, even if the system is very simple. Therefore, for such non-local partial differential equations, we only consider numerical solutions.
\begin{figure}
\begin{minipage}[!h]{0.49 \textwidth}
\centerline{\includegraphics[width=6cm,height=5cm]{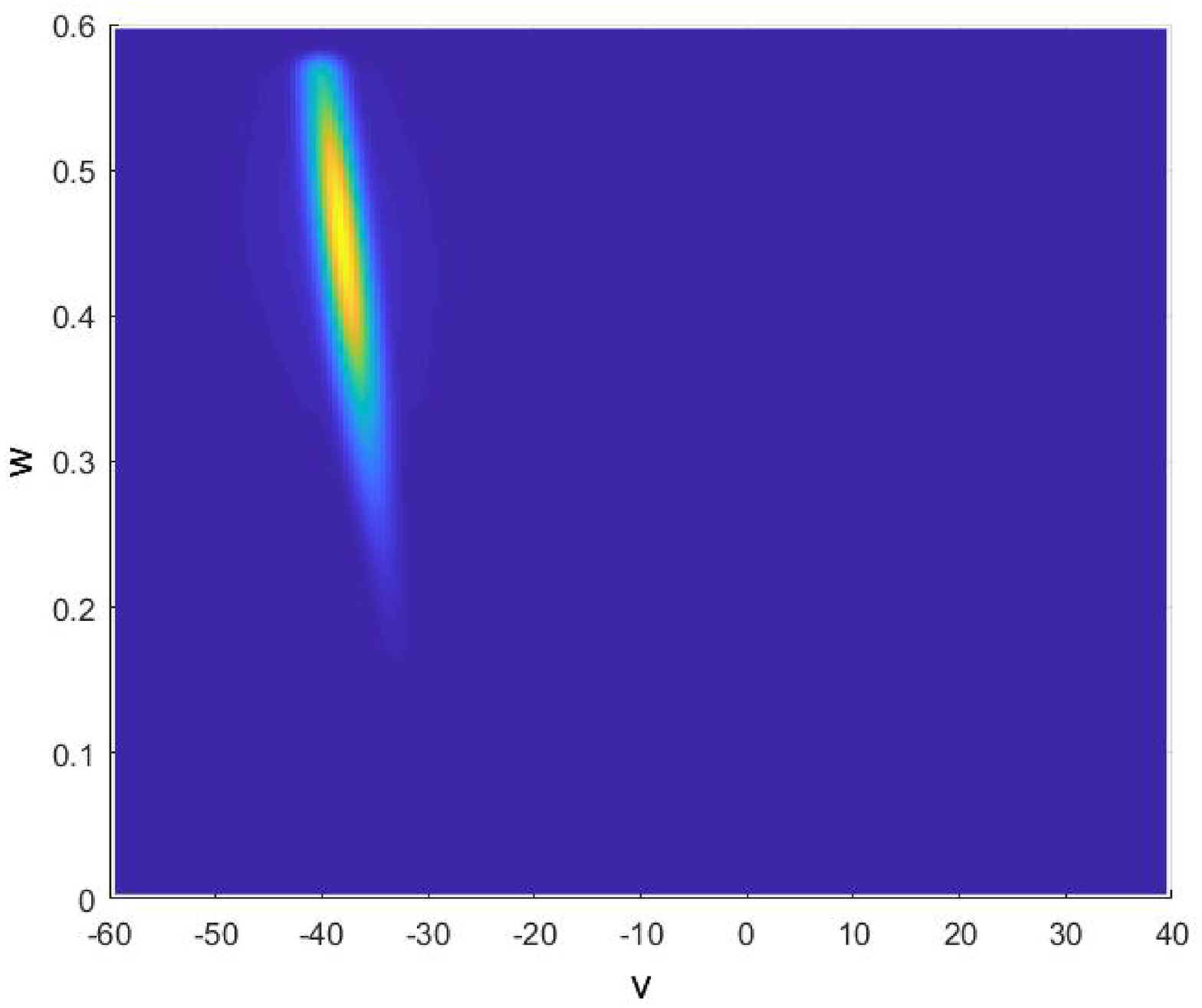}}
\centerline{(a) $t=1$}
\end{minipage}
\hfill
\begin{minipage}[!h]{0.49 \textwidth}
\centerline{\includegraphics[width=6cm,height=5cm]{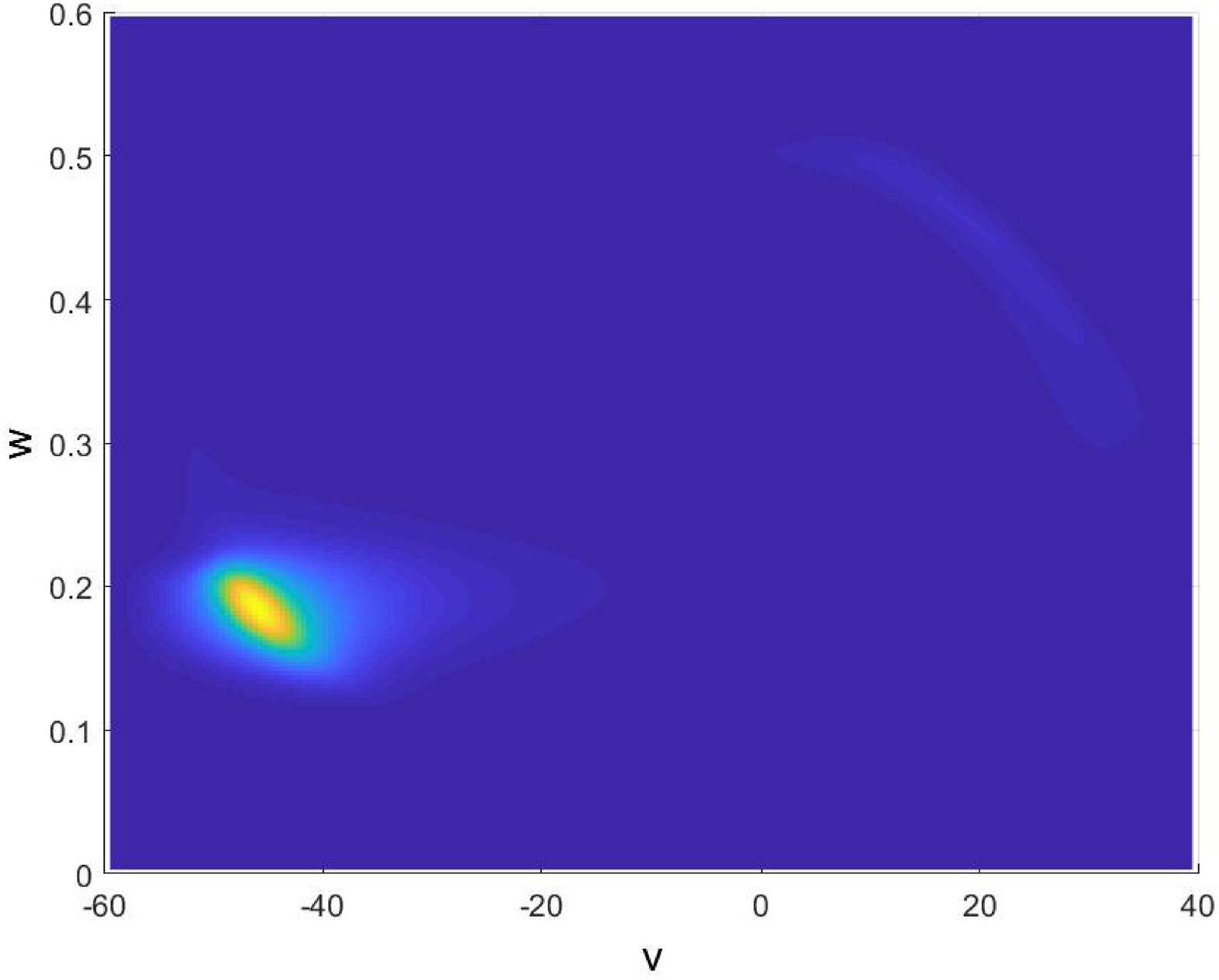}}
\centerline{(b) $t=20$}
\end{minipage}

\begin{minipage}[!h]{0.49 \textwidth}
\centerline{\includegraphics[width=6cm,height=5cm]{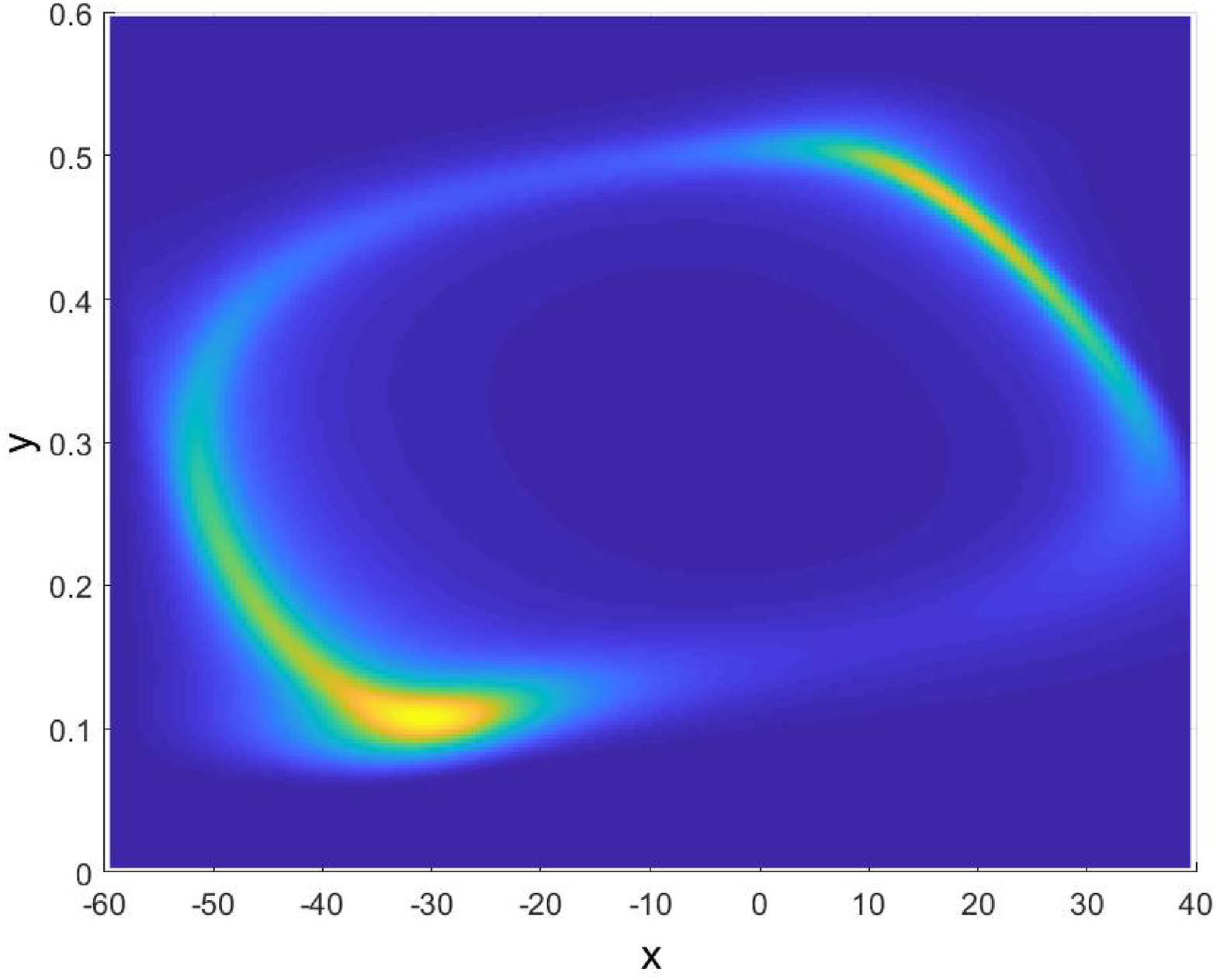}}
\centerline{(c) $t=70$}
\end{minipage}
\hfill
\begin{minipage}[!h]{0.49 \textwidth}
\centerline{\includegraphics[width=6cm,height=5cm]{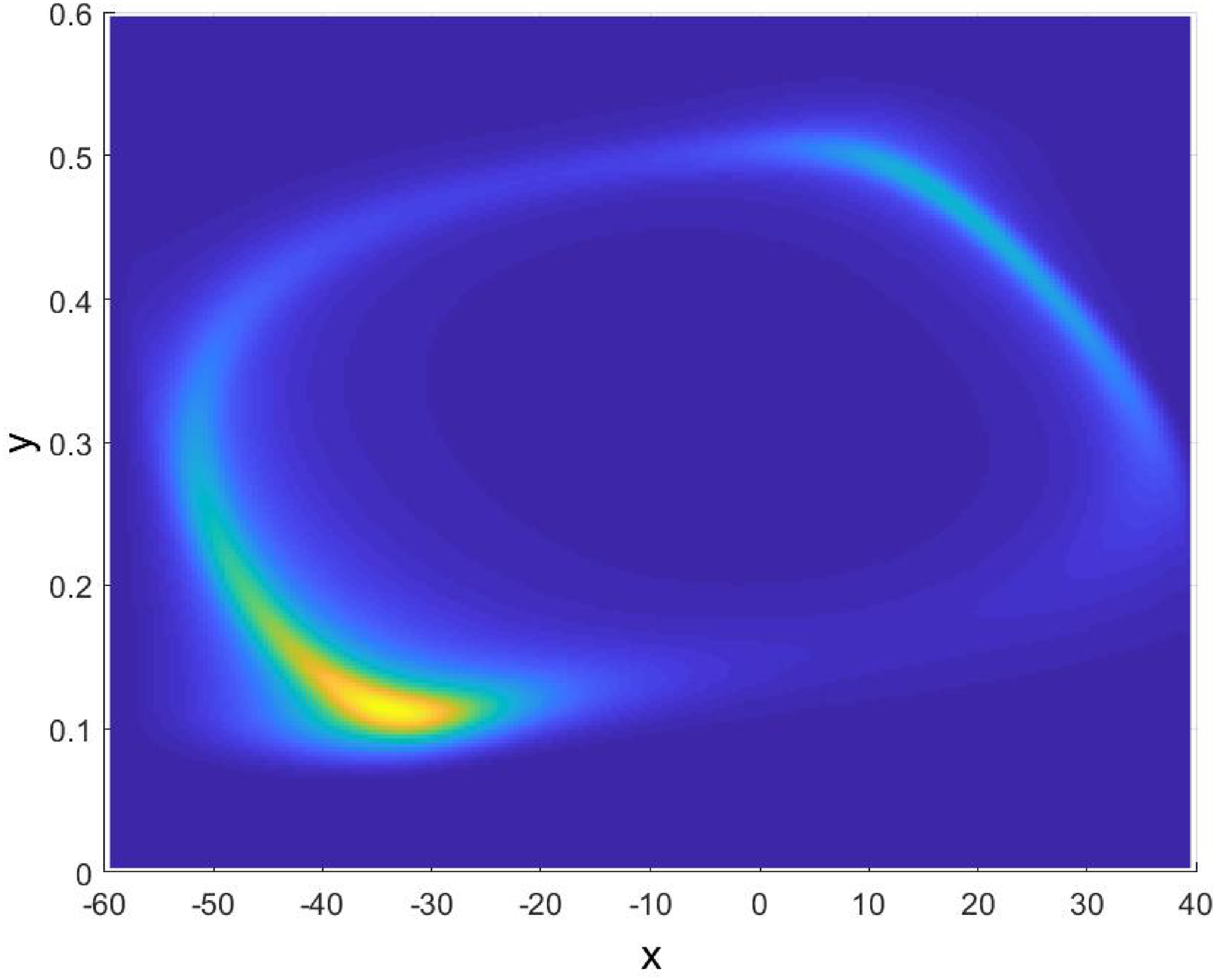}}
\centerline{(d) $t=100$}
\end{minipage}

\vfill
\caption{The probability density function of the stochastic Morris-Lecar model driven by L\'evy noise in the selected region at different time with L\'evy motion index $\alpha=0.5$ and noise intensity $\sigma=0.25$. The initial point is $(-32.7,0.4578)$. (a)-(d) shows the representative probability density function images of $t=1,~t=20,~t=70$ and $t=100$, respectively.}
\label{fig:2}
\end{figure}

In our case, the region $D = (-60, 40)\times(0, 0.6)$ is chosen to compute the probability density function.
Figure \ref{fig:2} shows the probability density function graph of the trajectory starting from the point $(-32.7,0.4578)$ at four fixed time. The initial point is on the big stable limit cycle and the L\'evy motion index $\alpha=0.5$ and noise intensity $\sigma=0.25$. The dark red region in Figure \ref{fig:2} represents the maximum value of the probability density function at that moment. At time $t=1$, the maximum value of the probability density function still gathers near the initial point. Then at $t=20$, the maximum value of the probability density function moves along the limit cycle. The maximum value along the trajectory of the limit cycle is more clearly revealed from the figure of $t=70$ and $t=100$. Obviously, the red area goes along the big stable limit cycle as time goes by. It shows that under this L\'evy motion index and noise intensity, the solution trajectory from this initial point does not change its state, and the system continues to maintain sustained oscillation. Only during these four fixed times, we can't judge whether the L\'evy noise in this case can cause the system to transfer state by the deadline of $t=100$. What we need is to observe the change in the maximum value of the probability density function in continuous time, that is the maximal likely trajectory.

\section{Result}\label{R}

\begin{figure}
\begin{minipage}[!h]{0.49 \textwidth}
\centerline{\includegraphics[width=6cm,height=5cm]{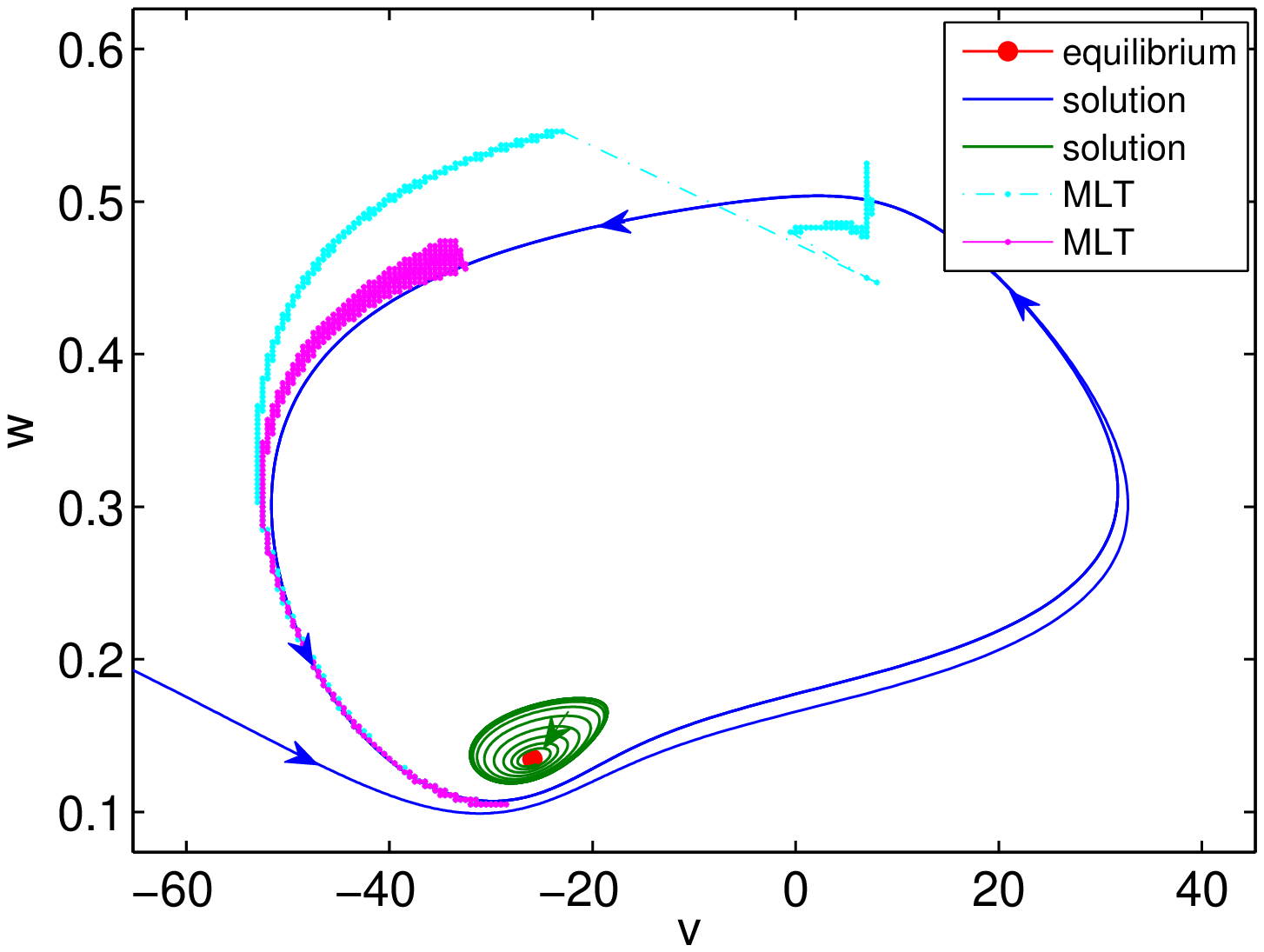}}
\centerline{(a) $\alpha=0.5, \sigma=0.25$}
\end{minipage}
\hfill
\begin{minipage}[!h]{0.49 \textwidth}
\centerline{\includegraphics[width=6cm,height=5cm]{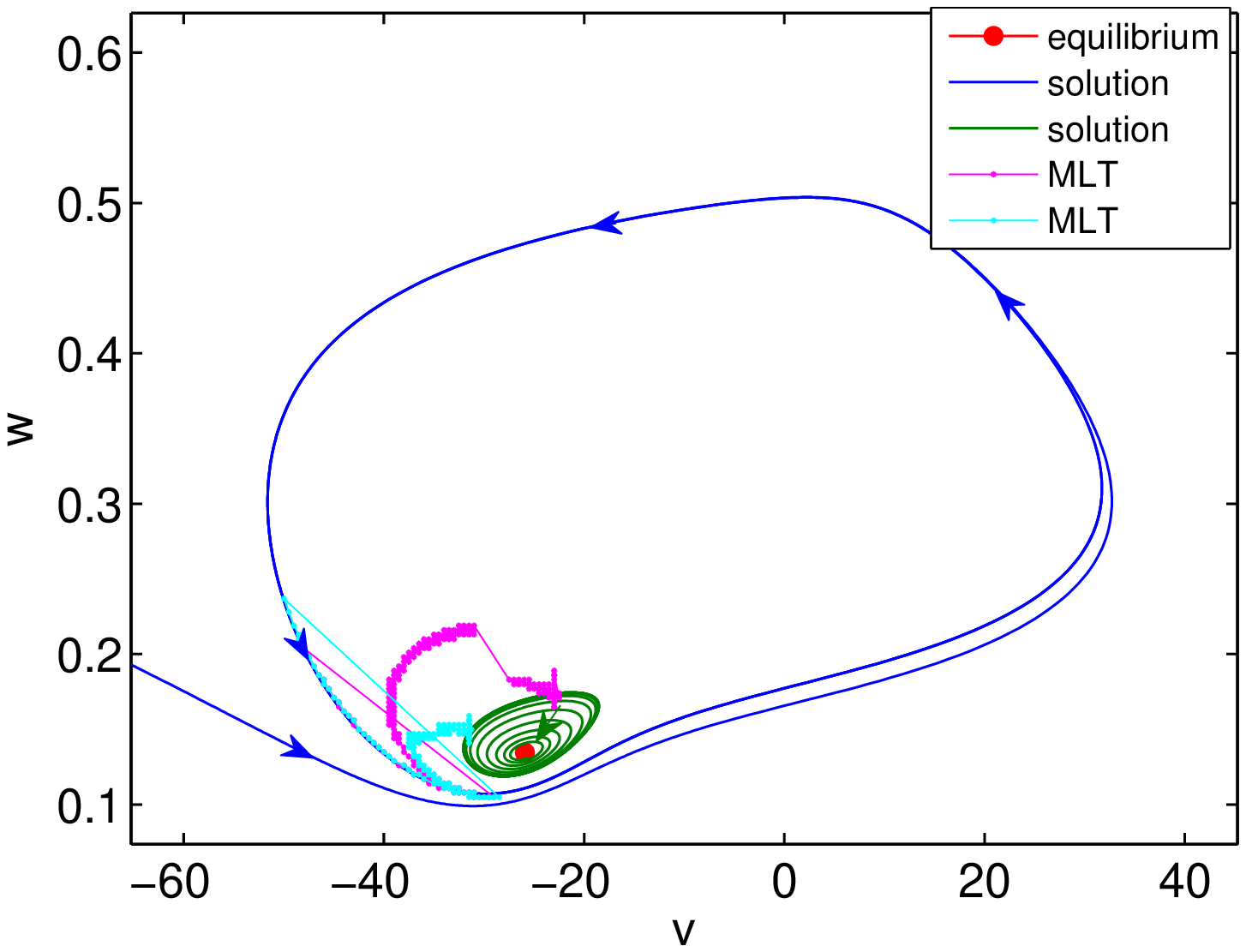}}
\centerline{(b) $\alpha=0.5, \sigma=0.25$}
\end{minipage}

\begin{minipage}[!h]{0.49 \textwidth}
\centerline{\includegraphics[width=6cm,height=5cm]{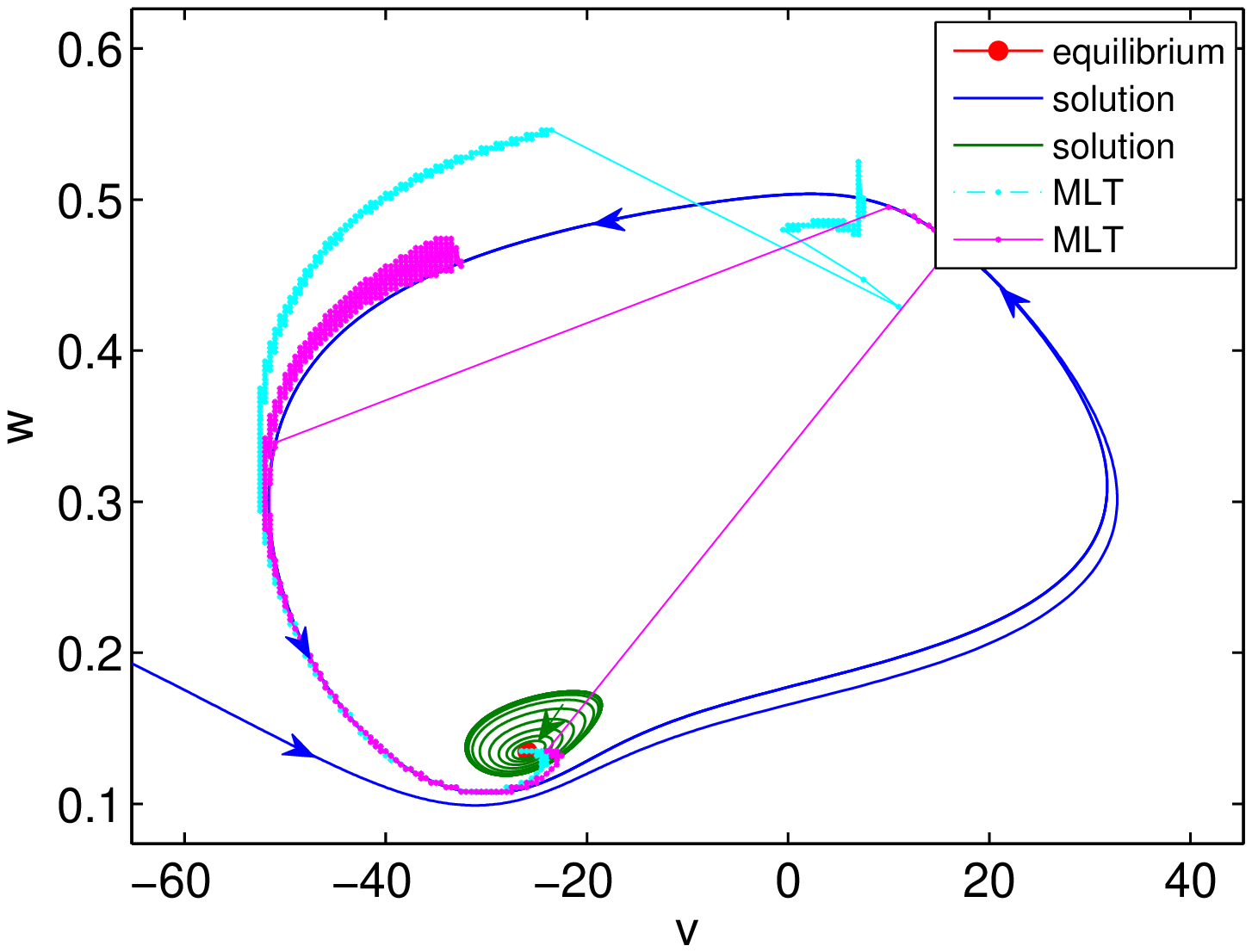}}
\centerline{(c) $\alpha=1.5, \sigma=0.25$}
\end{minipage}
\hfill
\begin{minipage}[!h]{0.49\textwidth}
\centerline{\includegraphics[width=6cm,height=5cm]{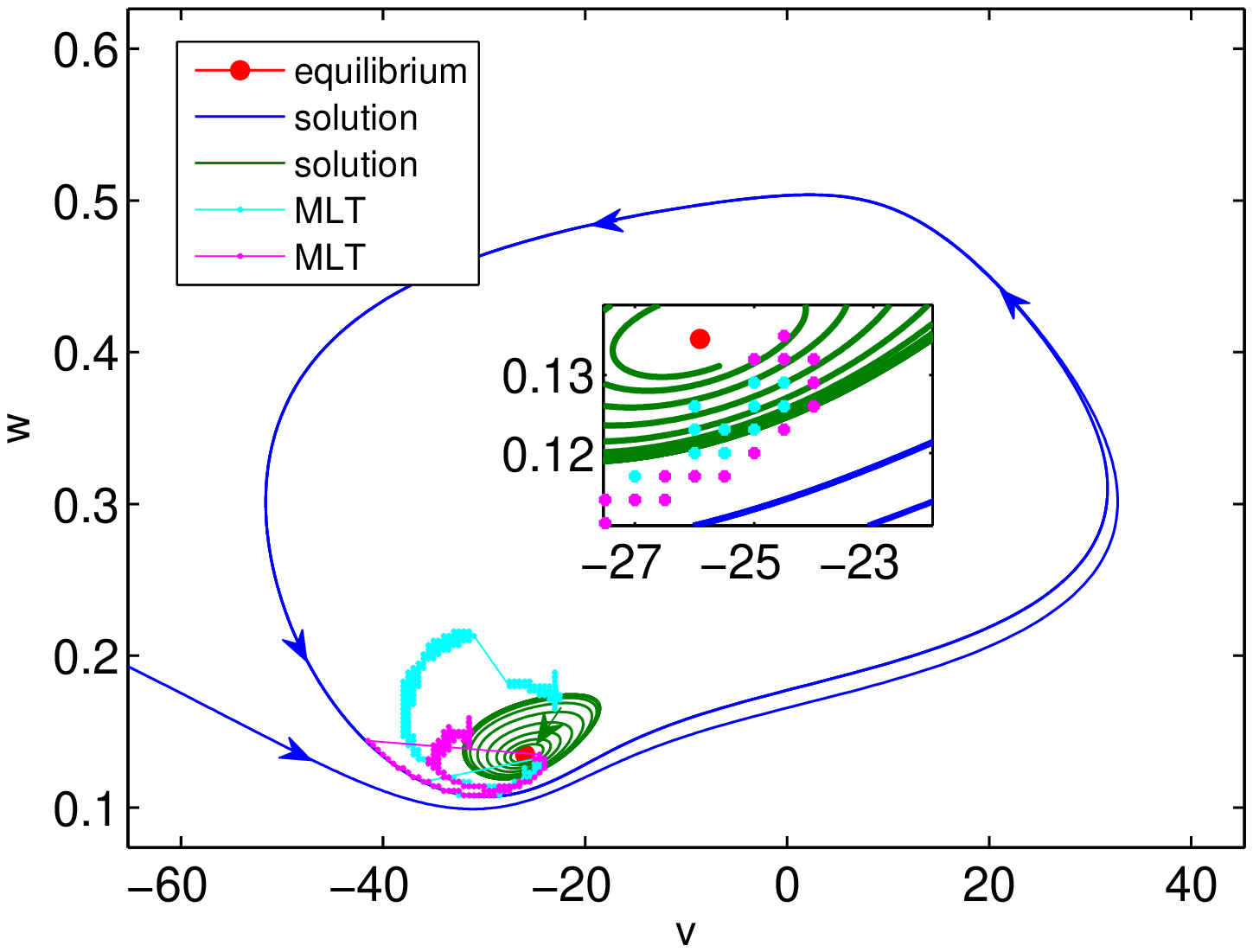}}
\centerline{(d) $\alpha=1.5, \sigma=0.25$}
\end{minipage}

\vfill
\caption{The maximal likely trajectories starting from different initial points in two dimensional plane. (a)The two initial points (-32.7,0.4578) and (7.459,0.5004) are on the stable limit cycle, the stochastic Morris-Lecar model with L\'evy motion index $\alpha=0.5$ and noise intensity $\sigma=0.25$.  (b)The two initial points (-22.73, 0.174) and (-31.27,0.15) are on the unstable limit cycles, the stochastic Morris-Lecar model with L\'evy motion index $\alpha=0.5$ and noise intensity $\sigma=0.25$. (c)The two initial points (-32.7,0.4578) and (7.459,0.5004) are on the stable limit cycle, the stochastic Morris-Lecar model with L\'evy motion index $\alpha=1.5$ and noise intensity $\sigma=0.25$. (d) The two initial points (-22.73, 0.174) and (-31.27,0.15) are on the unstable limit cycles, the stochastic Morris-Lecar model with L\'evy motion index $\alpha=1.5$ and noise intensity $\sigma=0.25$. The deadline is $t=100$.}
\label{fig:3}
\end{figure}

Figure \ref{fig:3} presents several maximal likely trajectories starting from four initial points with different L\'evy motion index and same noise intensity. Two initial points are on the stable limit cycle and the other two initial points are on the unstable limit cycle. For convenience of comparison, there are two kinds of L\'evy motion index with same noise intensity $\sigma=0.25$, one is $\alpha=0.5$ , the other is $\alpha=1.5$. The following summarizes the results of the comparison:

1, By comparing Figures \ref{fig:3} (a) with (b), we can see that under the same noise condition, the state of the two solution trajectories starting from the stable limit cycle has not changed, the maximal likely trajectories are still on the stable limit cycle in the time interval [0,100]. At the same time, the stable limit cycle attracts the solution trajectories starting from the boundary. That is to say, under the interference of such noise conditions ($\alpha=0.5, \sigma=0.25$), the system continues oscillating without state transition. While the borderline state of the system becomes the sustained oscillating state.

2, By comparing Figures \ref{fig:3} (c) with (d), the opposite conclusion can be obtained to the above. The state of the two solution trajectories starting from the stable limit cycle has changed. And the two maximal likely trajectories enter the attraction basin of the stable fixed point. This indicates that the noise under this condition induces the Morris-Lecar system to occur a state transition. Similarly, the fixed point attracts the two maximal likely trajectories starting from the unstable limit cycle. That is to say, under the interference of such noise conditions ($\alpha=1.5, \sigma=0.25$), the system changes from the sustained oscillating state to the resting state, and a state transition occurs. Similarly, the borderline state of the system also changes to the resting state under the disturbance of noise.

3, By comparing Figures \ref{fig:3} (a) with (c), it is clear to find that under the same initial point and noise intensity conditions, different L\'evy motion index lead to completely different phenomenon of the two maximal likely trajectories. When $\alpha=0.5$ ($\leq1$) the noise does not change the state behavior of the two maximal likely trajectories (the system remains in sustained oscillating state), but when $\alpha=1.5$ ($\geq1$), the two maximal likely trajectories change the oscillating state under the interference of L\'evy noise and enter the resting state.

4, By comparing Figures \ref{fig:3} (b) with (d), under the same initial point and noise intensity conditions, the two maximal likely trajectories starting from the small unstable limit cycle are similar to the case 3 and are attracted by the stable limit cycle when $\alpha=0.5$ ($\leq1$). That means the system changes from the borderline state to the oscillating state. While when $\alpha=1.5$ ($\geq1$), the two maximal likely trajectories are attracted by the stable point and the system becomes the resting state.

\begin{figure}[ht]
\centering %使插入的图片居中显示
\includegraphics[height=6cm ,width=8cm]{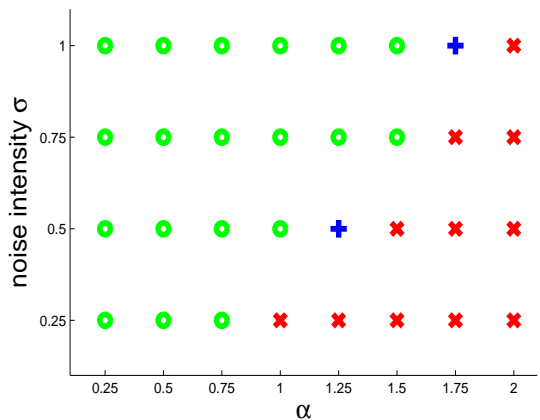}
\caption{Influence of L\'evy motion index $\alpha$ and noise intensity $\sigma$ on the state transition of two maximal likely trajectories. The initial points are (-32.7,0.4578) and (7.459,0.5004) on the stable limit cycle. If the two maximal likely trajectories do not enter the basin of the stable point (equilibrium) at the corresponding L\'evy motion index and noise intensity, we mark the point with a green circle ``$o$". If the two maximal likely trajectories enter the basin of the stable point (equilibrium), we mark the point with red ``x". If one maximal likely trajectory enters the basin of the stable point (equilibrium), the other does not, we mark the point blue ``$+$". As a comparison with the L\'evy situation, the independent Brownian motion $B_{t}$ is considered (corresponds to $\alpha=2$). The deadline is $t=100$.} % 插入图片的标题，一般放在图片的下方，放在表格的上方
\label{fig:4}
\end{figure}

\begin{figure}[ht]
\centering %使插入的图片居中显示
\includegraphics[height=6cm ,width=8cm]{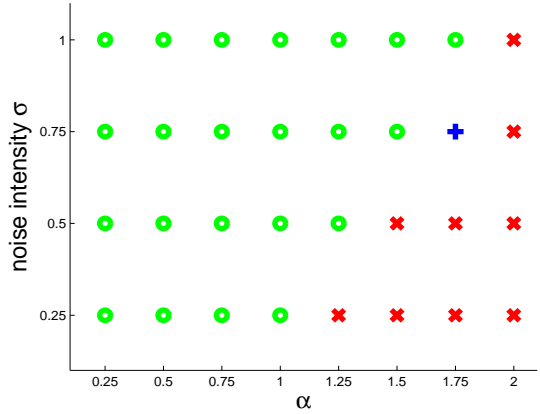}
\caption{Influence of L\'evy motion index $\alpha$ and noise intensity $\sigma$ on the state transition of two maximal likely trajectories. The initial points are (-22.73,0.174) and (-31.27,0.15) on the small unstable limit cycle. If the two maximal likely trajectories do not enter the basin of the stable point (equilibrium) at the corresponding L\'evy motion index and noise intensity (which means the two maximal likely trajectories attracted by the big stable limit cycle), we mark the point with a green circle ``$o$". If the two maximal likely trajectories enter the basin of the stable point (equilibrium), we mark the point with red ``x". Similarly, if one maximal likely trajectory enters the basin of the stable point (equilibrium), the other does not, we mark the point blue ``$+$". As a comparison with the L\'evy situation, the independent Brownian motion $B_{t}$ is considered (corresponds to $\alpha=2$).The deadline is $t=100$.} % 插入图片的标题，一般放在图片的下方，放在表格的上方
\label{fig:5}
\end{figure}

From what has been discussed above, two issues are worth paying attention to: under what noise conditions, the system will change from the oscillating state to the resting state. Which is more attractive for the most likely trajectory starting from the borderline state, the oscillating state or the resting state? Figure \ref{fig:4} and Figure \ref{fig:5} answer the two questions.

Figure \ref{fig:4} and Figure \ref{fig:5} present a summary of the behaviors of the four maximal likely trajectories under different L\'evy motion index and noise intensities. The difference between Figure \ref{fig:4} and Figure \ref{fig:5} is that the initial point is different. The starting point of Figure \ref{fig:4} is on the big stable limit cycle, this means that the Morris-Lecar system is in the oscillating state. And the corresponding starting point of Figure \ref{fig:5} is on the small unstable limit cycle, this means the Morris-Lecar system is in the borderline state.

From Figure \ref{fig:4}, it can be seen that when $\alpha=1$, the state of the maximal likely trajectories starting from the stable limit cycle begins to change (the critical case is $\alpha=1$). When the noise intensity $\sigma=0.25$, the sustained oscillating state is changed by the interference of noise when $\alpha \geq 1$. As the noise intensity increases, the L\'evy noise with larger $\alpha$ can change the sustained oscillating state into the resting state. Similar phenomena also occur in the two maximal likely trajectories starting from the unstable limit cycle. From Figure \ref{fig:5}, unlike the other two maximal likely trajectories starting on the big stable limit cycle, the two maximal likely trajectories are attracted by the stable point when $\alpha=1.25$, so the turning point is $\alpha=1.25$ . Similarly, when the noise intensity $\sigma=0.25$, the state of the two maximal likely trajectories eventually enter the resting state when $\alpha \geq 1.25$. As the noise intensity increases, the larger $\alpha$ is needed to make the borderline state change to the resting state.

It is worth noting that the Brownian case ($\alpha=2$) presents a phenomenon that is completely different from the non-Gaussian L\'evy case. Under the four selected noise intensities conditions, the four maximal likely trajectories are all attracted by the stable points and this system enter the resting state.

It is known that abnormal discharges in the brain may cause seizures and hypoxia in the brain, resulting in a series of adverse reactions. Sustained oscillation may correspond to abnormal discharge behavior. In this paper, we consider the effect of L\'evy noise on the discharge behavior of neuronal system (transition between two states), this is of great significance for us to fully understand some characteristics of isolated neurons.

\section{Conclusion}\label{CON}

We have analyzed the state transition of the Morris-Lecar neural model perturbed by a symmetric $\alpha$-stable L\'evy noise (non-Gaussian noise). A new perspective has been provided to describe two state transitions: one is from the oscillating state to the resting state and the other one is state selection of trajectories starting from the borderline state, by means of the maximal likely trajectory. The effects of the L\'evy motion index $\alpha$ and the noise intensity $\sigma$ on the state of the the maximal likely trajectories have been presented by numerical simulations. We found that smaller jumps of the L\'evy motion and relatively smaller noise intensity are conducive to the state transition of the Morris-Lecar neural model from the oscillating state to the resting state, while higher noise intensity and larger jumps of the L\'evy motion promote the system from the borderline state change to the oscillating state.

Moreover, Brownian motion has been considered in comparison. Compared with the L\'evy case, the continuous Brownian motion showed a great difference when the noise intensity is fixed: under the four fixed noise intensity we select, the system disturbed by Brownian noise enter into the resting state whether it is in the state of oscillation or in the state of boundary. In fact, the maximal likely trajectory approach allow us to capture the primary statistical features and geometric behaviors of the solution trajectories in the stochastic Morris-Lecar neural model.

\section{Acknowledgements}

We would like to thank Dr. Xiaoli Chen for helpful discussions. This work was partly supported
by the National Science Foundation Grant (NSF) No. 1620449, and the National Natural Science
Foundation of China (NSFC) Grant Nos. 11531006 and 11771449.

%%%%%%%%%%%%%%%%%%%%%%%%%%%%%%%%%%%%%%%%%%%‘’‘’‘’‘’‘’‘’‘’‘’‘’‘’‘’‘’‘’‘’‘’‘’‘’‘’‘’‘’‘’‘’‘’‘’‘’‘’‘’‘’‘’‘’‘’‘’‘’‘’‘’‘’‘’‘’‘’‘’‘’‘’‘’‘’‘’‘’‘’‘’‘’‘’‘’‘’‘’‘’‘’‘’‘’‘’‘’‘’‘’‘’‘’‘’‘’‘’‘’‘’‘’‘’‘’‘’‘’‘’‘’‘’‘’‘’‘’‘’‘’‘’‘’‘’‘’‘’‘’‘’‘’‘’‘’‘’‘’‘’‘’‘’‘’‘’‘’‘’‘’‘’‘’‘’‘’‘’‘’‘’‘’‘’‘’‘’‘’‘’‘’‘’‘’‘’‘’‘’‘’‘’‘’‘’‘’‘’‘’‘’‘’‘’‘’‘’‘’‘’‘’‘’‘’‘’‘’‘’‘’‘’‘’‘’‘’‘’‘’‘’‘’‘’‘’‘’‘’‘’‘’‘’‘’‘’‘’‘’‘’‘’‘’‘’‘’‘’‘’‘’‘’‘’‘’‘’‘’‘’‘’‘’‘’‘’‘’‘’‘’‘’‘’‘’‘’‘’‘’‘’‘’‘’‘’‘’‘’‘’‘’‘’‘’‘’‘’‘’‘’‘’‘’‘’‘’‘’‘’‘’‘’‘’‘’‘’‘’‘’‘’‘’‘’‘’‘’‘’‘’‘’‘’‘’‘’‘’‘’‘’‘’‘’‘’‘’‘’‘’‘’‘’‘’‘’‘’‘’‘’‘’‘’‘’‘’‘’‘’‘’‘’‘’‘’‘’‘’‘’‘’‘’‘’‘’‘’‘’‘’‘’‘’‘’‘’‘’‘’‘’‘’‘’‘’‘’‘’‘’‘’‘’‘’‘’‘’‘’‘’‘’‘’‘’‘’‘’‘’‘’‘’‘’‘’‘’‘’‘’‘’‘’‘’‘’‘’‘’‘’‘’‘’‘’‘’‘’‘’‘’‘’‘’‘’‘’‘’‘’‘’‘’‘’‘’‘’‘’‘’‘’‘’‘’‘’‘’‘’‘’‘’‘’‘’‘’‘’‘’‘’‘’‘’‘’‘’‘’‘’‘’‘’‘’‘’‘’‘’‘’‘’‘’‘’‘’‘’‘’‘’‘’‘’‘’‘’‘’‘’‘’‘’‘’‘’‘’‘’‘’‘’‘’‘’‘’‘’‘’‘’‘’‘’‘’‘’‘’‘’‘’‘’‘’‘’‘’‘’‘’‘’‘’‘’‘’‘’‘’‘’‘’‘’‘’‘’‘’‘’‘’‘’‘’‘’‘’‘’‘’‘’‘’‘’‘’‘’‘’‘’‘’‘’‘’‘’‘’‘’‘’‘’‘’‘’‘’‘’‘’‘’‘’‘’‘’‘’‘’‘’‘’‘’‘’‘’‘’‘’‘’‘’‘

\bibliographystyle{unsrt}
\bibliography{levy}

\end{document}